\documentstyle[12pt]{article}
\def\Bbb R{{\rm \bf R}}
\def\proclaim#1{\vskip2mm{\bf #1}\em}
\def\endproclaim{\em \vskip2mm}
\def\tag#1{\eqno(#1)}
\def\gathered{\begin{array}{c}}
\def\endgathered{\end{array}}
\def\text{\mbox}

\begin{document}

\title {Two Sides of Probe Method\\
and\\
Obstacle with Impedance Boundary Condition}
\author{Masaru IKEHATA\\
Department of Mathematics,
Faculty of Engineering\\
Gunma University, Kiryu 376-8515, JAPAN\\
\small{(Dedicated to Professor Minoru Murata on the occasion of his 60th birthday)}
}
\date{}
\maketitle
\begin{abstract}
An inverse boundary value problem for the Helmholtz equation in a
bounded domain is considered.  The problem is to extract
information about an unknown obstacle embedded in the domain with
unknown impedance boundary condition (the Robin condition) from
the associated Dirichlet-to-Neumann map. The main result is a
characterization of the unknown obstacle via the sequences that
are constructed by the Dirichlet-to-Neumann map, under smallness
conditions on the wave number and the upper bound of the
impedance. Moreover two alternative simple proofs of a previous result of
Cheng-Liu-Nakamura which are based on only some energy estimates,
an analysis of the blowup of the energy of so-called reflected solutions
and an application of the enclosure method to the problem are also
given.

\noindent
AMS: 35R30

\end{abstract}

\section{Introduction}

In this paper, we consider an inverse boundary value problem for the Helmholtz equation in a bounded domain.
The problem is to extract information about an unknown obstacle
embedded in the domain with unknown impedance boundary condition from
the associated Dirichlet-to-Neumann map.

Let $\Omega$ be a bounded domain in $\Bbb R^m(m=2,3)$ with Lipschitz boundary.
Let $D$ be an open subset with Lipschitz boundary of $\Omega$ and
satisfy that $\overline D\subset\Omega$; $\Omega\setminus\overline D$ is connected.
We always assume that $D$ is given by a union of finitely many bounded Lipschitz domains
$D_1,\cdots, D_N$ such that $\overline D_j\cap\overline D_l=\emptyset$ if $j\not=l$.

\noindent
We denote by $\nu$ the unit outward normal relative to $D$.
Let $k\ge 0$.  We always assume that $0$ is not a Dirichlet eigenvalue of $\triangle +k^2$ in $\Omega$.
Let $\lambda\in L^{\infty}(\partial D)$.

Given $f\in H^{1/2}(\partial\Omega)$ we say that $u\in H^1(\Omega\setminus\overline D)$ is
a weak solution of the elliptic problem
$$\begin{array}{c}
\displaystyle
\triangle u+k^2u=0\,\,\text{in}\,\Omega\setminus\overline D,\\
\\
\displaystyle
\frac{\partial u}{\partial\nu}+\lambda u=0\,\,\text{on}\,\partial D,\\
\\
\displaystyle
u=f\,\,\text{on}\,\partial\Omega
\end{array}
\tag {1.1}
$$
if $u$ satisfies $u=f$ on $\partial\Omega$ in the sense of trace
and, for all $\varphi\in H^1(\Omega\setminus\overline D)$ with $\varphi=0$ on $\partial\Omega$
in the sense of trace
$$\displaystyle
\int_{\partial D}\lambda u\varphi dS-\int_{\Omega\setminus\overline D}\nabla u\cdot\nabla\varphi dy
+k^2\int_{\Omega\setminus\overline D}u\varphi dy=0.
$$
We assume that: there exists a positive constant $C>0$ such that
$\text{Im}\,\lambda(x)\ge C$ for almost all $x\in\partial D$. This
assumption is motivated by a possibility of application to inverse
scattering problems of electromagnetic/acoustic
wave(\cite{CK,CP,GK}). It is really routine to see that the weak solution of
(1.1) exists and is unique.

Define the bounded linear functional $\Lambda_Df$ on $H^{1/2}(\partial\Omega)$
by the formula
$$\displaystyle
<\Lambda_Df,h>
=-\int_{\partial D}\lambda u vdS
+\int_{\Omega\setminus\overline D}\nabla u\cdot\nabla vdy-
k^2\int_{\Omega\setminus\overline D}uv dy,\,\,h\in H^{1/2}(\partial\Omega)
$$
where $u$ is the weak solution of (1.1) and $v\in H^1(\Omega\setminus\overline D)$ is an arbitrary function with
$v=h$ on $\partial\Omega$ in the sense of the trace.
The map $\Lambda_D:f\longmapsto \Lambda_Df$ is called the Dirichlet-to-Neumann map.
We set $\Lambda_D=\Lambda_0$ in the case when $D=\emptyset$.
\noindent

\noindent
Here we consider the problem

{\it\noindent Inverse Problem}
Extract information about the shape and location of $D$
from $\Lambda_D$ or its partial knowledge.

In \cite{I1} the author introduced the probe method which gives a
general idea to obtain a reconstruction formula of unknown
objects embedded in a known background medium from a mathematical
counterpart(the Dirichlet-to-Neumann map) of the measured data of
some physical quantity on the boundary of the medium. The method
was applied to an inverse boundary value problem in elasticity
\cite{I4} and mixed problems \cite{DIN}.

In \cite{I3} we considered {\it Inverse Problem} in the two
extreme cases when $\lambda\longrightarrow
0$(sound-hard)/$\lambda\longrightarrow\infty$(sound-soft). Using
the probe method, we established reconstruction formulae of
$\partial D$ by using $(\Lambda_0-\Lambda_D)f$ for infinitely many
$f$ provided both $\partial\Omega$ and $\partial D$ are $C^2$. In
\cite{I2} we considered inverse obstacle scattering problems with
a fixed frequency and established a way of reconstructing
sound-hard/sound-soft obstacle from the scattering data.

Recently Erhard-Potthast (\cite{EP}) studied the probe method
numerically. They considered
{\it Inverse Problem}
in the case when $\lambda\longrightarrow \infty$ and
computed an approximation of the corresponding {\it indicator function} by employing the
techniques of the point source and singular sources methods by Potthast (\cite{P1,P2}).
Applying the probe method to
{\it Inverse Problem}
for general $\lambda$, Cheng-Liu-Nakamura \cite{CLN}
gave a corresponding reconstruction formula of $\partial D$ provided also both $\partial\Omega$ and $\partial D$
are $C^2$ and $\lambda$ is unknown.  Quite recently a numerical testing
of the probe method applied to {\it Inverse Problem} has been done in \cite{CLN2}.

By the way, in order to give an explanation for an observation
obtained in \cite{EP}, quite recently the author gave a new
formulation of the probe method, raised new questions about the
probe method itself and gave answers to some of them in the two
extreme cases \cite{I6}.

The aim of this paper is to: reconsider the application of the previous version of the probe method to
{\it Inverse Problem} and add a new result based on the new formulation of the probe method
in \cite{I6}.
We clarify that the probe method has two different sides.  One side has a common character
with the enclosure method which was introduced in \cite{IE} and is based on the blowup
property of the exponentially growing solutions.  Thus we apply also the enclosure method
to {\it Inverse Problem} in the final section for comparison.

For application of the previous version of the probe method to inverse boundary value problems related
to cracks (non volumetric discontinuity) see \cite{IN, NUW}.

\section{Two Sides of Probe Method}

In \cite{I6} we introduced a simpler formulation of the probe
method and further investigated the method itself in the two
extreme cases when $\lambda\longrightarrow 0$ or
$\lambda\longrightarrow\infty$. We start with describing the
formulation.

Given a point $x\in\Omega$ let $N_x$ denote the set of all piecewise linear
curves $\sigma:[0,\,1]\longmapsto\overline\Omega$ such
that :

$\bullet$  $\sigma(0)\in\partial\Omega$, $\sigma(1)=x$ and
$\sigma(t)\in\Omega$ for all $t\in]0,\,1[$;

$\bullet$  $\sigma$ is injective.

\noindent
We call $\sigma\in N_x$ a {\it needle with tip} at $x$.

\noindent
We denote by $G_k(x)$ the standard fundamental solution of the Helmholtz equation.
For the new formulation of the probe method we need the following.

{\bf\noindent Definition 2.1.}
Let $\sigma\in N_x$.  We call the sequence $\xi=\{v_n\}$
of $H^1(\Omega)$ solutions of the Helmholtz equation
a {\it needle sequence} for $(x,\sigma)$ if it satisfies,
for each
fixed compact set $K$ of $\Bbb R^m$ with $K\subset\Omega\setminus
\sigma([0,\,1])$
$$\displaystyle
\lim_{n\longrightarrow\infty}(\Vert v_n(\,\cdot\,)-G_k(\,\cdot\,-x)\Vert_{L^2(K)}
+\Vert\nabla\{v_n(\,\cdot\,)-G_k(\,\cdot\,-x)\}\Vert_{L^2(K)})=0.
$$
\noindent
The existence of the needle sequence has been ensured in \cite{I3}.

{\bf\noindent 2.1. Side A of Probe Method-A Review on Cheng-Liu-Nakamura's Result}

\noindent
{\bf\noindent Definition 2.2.}
Given $x\in\,\Omega$, needle $\sigma$ with tip $x$
and needle sequence $\xi=\{v_n\}$ for $(x,\sigma)$
define
$$\displaystyle
I(x,\sigma,\xi)_n=\text{Re}\,<(\Lambda_0-\Lambda_D)f_n,\overline f_n>,\,n=1,2,\cdots
$$
where
$$\displaystyle
f_n(y)=v_n(y),\,\,y\in\,\partial\Omega.
$$

\noindent $\{I(x,\sigma,\xi)_n\}_{n=1,2,\cdots}$ is a sequence
depending on $\xi$ and $\sigma\in N_x$.
We call the sequence the {\it indicator sequence}.

Now we can say that the {\it probe method} is a method of {\it probing} inside
a given medium by monitoring the behaviour of the indicator sequences
for many needles.  For the description of the behaviour we introduce a function
defined outside the obstacle.

{\bf\noindent Definition 2.3.}
The {\it indicator function} $I$ is defined by the formula
$$\begin{array}{c}
\displaystyle
I(x)
=\int_{\Omega\setminus\overline D}\vert\nabla w_x\vert^2dy
-k^2\int_{\Omega\setminus\overline D}\vert w_x\vert^2 dy\\
\\
\displaystyle
+\int_{D}\vert\nabla G_k(y-x)\vert^2 dy-k^2\int_{D}\vert G_k(y-x)\vert^2 dy\\
\\
\displaystyle
+\int_{\partial D}\{(\text{Re}\,\lambda)(\vert G_k(\,\cdot\,-x)\vert^2-\vert w_x\vert^2)
-2(\text{Im\,}\lambda)\text{Im}\,\{w_x\overline G_k(\,\cdot\,-x)\}\}dS,\,\,x\in\Omega\setminus\overline D
\end{array}
\tag {2.1}
$$
where $w_x\in H^1(\Omega\setminus\overline D)$ is the unique weak solution of the problem
$$\begin{array}{c}
\displaystyle
\triangle w+k^2w=0\,\,\text{in}\,\Omega\setminus\overline D,\\
\\
\displaystyle
\frac{\partial w}{\partial\nu}+\lambda w=-(\frac{\partial}{\partial\nu}G_k(\,\cdot\,-x)+\lambda
G_k(\,\cdot\,-x))
\,\,\text{on}\,\partial D,\\
\\
\displaystyle
w=0\,\,\text{on}\,\partial\Omega.
\end{array}
$$
The function $w_x$ is called the {\it reflected solution} by $D$.

The following theorem extends the result in \cite{I3} and is
nothing but a result established in \cite{CLN}.

\proclaim{\noindent Theorem A.}
We have:

$\bullet$ (A.1) given $x\in\Omega\setminus\overline D$ and needle $\sigma$
with tip at $x$ if $\sigma(]0,1])\cap\overline D=\emptyset$,
then for any needle sequence $\xi=\{v_n\}$ for $(x,\sigma)$
the sequence $\{I(x,\sigma,\xi)_n\}$ converges to the indicator function $I(x)$;

$\bullet$ (A.2) for each $\epsilon>0$
$$\displaystyle
\sup_{\text{dist}\,(x,D)>\epsilon}I(x)<\infty;
$$

$\bullet$ (A.3) given point $a\in\partial D$
$$\displaystyle
\lim_{x\longrightarrow a}I(x)=\infty
$$
provided both $\partial\Omega$ and $\partial D$ are $C^2$ and $\lambda\in C^1(\partial D)$.

\endproclaim

\noindent
(A.2) is trivial. (A.3) is the most important property
of the indicator function. The indicator function blows up at the boundary
of the obstacle.  (A.1) gives a way of calculating the indicator function
by using the indicator sequence and is a direct consequence of the integral
identity (2.3) given below and the well posedness of the problem
$$\begin{array}{c}
\displaystyle
\triangle w+k^2w=0\,\,\text{in}\,\Omega\setminus\overline D,\\
\\
\displaystyle
\frac{\partial w}{\partial\nu}+\lambda w=g
\,\,\text{on}\,\partial D,\\
\\
\displaystyle
w=0\,\,\text{on}\,\partial\Omega
\end{array}
\tag {2.2}
$$
where $g\in H^{-1/2}(\partial D)$.

\proclaim{\noindent Proposition 2.1(\cite{CLN}).}
For all $f\in H^{1/2}(\partial\Omega)$
$$\begin{array}{c}
\displaystyle
<(\Lambda_0-\Lambda_{D})f,\overline f>
=\int_{\Omega\setminus\overline D}\vert\nabla(u-v)\vert^2dy
-k^2\int_{\Omega\setminus\overline D}\vert u-v\vert^2 dy\\
\\
\displaystyle
+\int_{D}\vert\nabla v\vert^2 dy-k^2\int_{D}\vert v\vert^2 dy\\
\\
\displaystyle
+\int_{\partial D}\{(\text{Re}\,\lambda)(\vert v\vert^2-\vert u-v\vert^2)
-2(\text{Im\,}\lambda)\text{Im}\,((u-v)\overline v)\}dS
+i\int_{\partial D}\text{Im}\,\lambda \vert u\vert^2 dS
\end{array}
\tag {2.3}
$$
where $u$ solves (1.1); $v$ solves
$$\begin{array}{c}
\displaystyle
(\triangle +k^2)v=0\,\,\text{in}\,\,\Omega,\\
\\
\displaystyle
v=f\,\,\text{on}\,\,\partial\Omega.
\end{array}
$$
\endproclaim

\noindent
Note that: this identity is a generalization of the corresponding one in the case
when $\lambda\longrightarrow 0$ which has been established in \cite{I3}.

Using (A.1) to (A.3), one can define another indicator function
which gives a previous formulation of the probe method, however, we do not
repeat to do it (see \cite{I3} for the previous formulation of the probe method).

The proof of (A.3) given by Cheng-Liu-Nakamura is not trivial and quite involved.
Here we describe the review on their proof of (A.3) by using the formulation given here
and points out also the problem.

They start with deriving the estimate from (2.1):
$$\begin{array}{c}
\displaystyle
I(x)\ge \int_D\vert\nabla G_k(y-x)\vert^2 dy-k^2\int_D\vert G_k(y-x)\vert^2dy
-k^2\int_{\Omega\setminus\overline D}\vert w_x\vert^2dy\\
\\
\displaystyle
-C_1\int_{\partial D}\vert G_k(y-x)\vert^2dS(y)
-C_2\int_{\partial D}\vert w_x\vert^2 dy
\end{array}
\tag {2.4}
$$
where both $C_1\ge 0$ and $C_2\ge 0$ are independent of $x$.
Note that both $C_1$ and $C_2$ vanish in the case when $\lambda\longrightarrow 0$
and (2.4) becomes the inequality (25) in \cite{I3}.

It is easy to see that the second term of (2.4) is bounded as $x\longrightarrow a\in\partial D$.
Using the completely same argument as done in \cite{I3} in the case when $\lambda\longrightarrow 0$
and the regularity assumptions on $\partial\Omega$, $\partial D$ and $\lambda$
(say $\lambda\in C^1(\partial D)$)
which ensure the $H^2(\Omega\setminus\overline D)$
regularity of solutions $p\in H^1(\Omega\setminus\overline D)$ of the problem
$$\begin{array}{c}
\displaystyle
\triangle p+k^2p\in L^2(\Omega\setminus\overline D),\\
\\
\displaystyle
\frac{\partial p}{\partial\nu}+\lambda p=0\,\,\text{on}\,\partial D,\\
\\
\displaystyle
p=0\,\,\text{on}\,\partial\Omega,
\end{array}
$$
one knows that the $L^2$-norm of the reflected solution $w_x$ over $\Omega\setminus\overline D$
is bounded as $x\longrightarrow a$(Theorem 3 of \cite{CLN}).  Note that, in their proof,
it is just assumed that $\lambda\in L^{\infty}(\partial D)$, however, we do not know
how to deduce the desired regularity for $p$ under the condition.

\noindent
This means that
one may propose: {\it prove (A.3) under the regularity condition
$\lambda\in L^{\infty}(\partial D)$.}

For the proof of (A.3) they showed that
$$\displaystyle
\int_{\partial D}\vert w_x\vert^2 dS(y)=O(\int_{\partial D}\vert G_k(y-x)\vert^2dS(y))
\tag {2.5}
$$
as $x\longrightarrow a$
(Theorem 4 of \cite{CLN}).  Since
$$\displaystyle
\lim_{x\longrightarrow a}\int_D\vert\nabla G_k(y-x)\vert^2 dy=\infty
$$
and
$$\displaystyle
\lim_{x\longrightarrow a}\frac{\displaystyle\int_{\partial D}\vert G_k(y-x)\vert^2dS(y)}
{\displaystyle\int_D\vert\nabla G_k(y-x)\vert^2 dy}=0,
$$
from (2.4) one obtains the desired conclusion.  The proof of (2.5) is done by
extracting the leading term of $w_x$ as $x\longrightarrow a$ that contributes
the blowup of $\Vert w_x\vert_{\partial D}\Vert_{L^2(\partial D)}^2$ as $x\longrightarrow a$
and needs a careful analysis of $w_x$ in a neighbourhood of $a$ (4 pages in the journal!).

We point out that (2.5) is unnecessary for the proof of (A.3).
In Section 3 we present two alternative simple proofs of (A.3).
The proofs do not make use of the local expression of the function $w_x$ in
a neighbourhood of $a\in\partial D$ and are based on
some energy estimates only.  This means that it maybe possible
to apply the method to more general cases (other inverse boundary value problems for elliptic equations
and systems of equations) without essential changes and difficulty.
Those belong to our future study.

{\bf\noindent 2.2. Side B of Probe Method}

\noindent The result in this subsection is new and, needless to
say, not covered in \cite{CLN}. It is related to another side of
the probe method.

In order to describe another side of the probe method we introduce
two positive constants appearing in two types of the Poincar\'e
inequalities (e.g., see \cite{SS,Z} and also \cite{I6}). One is
given in the following.

\proclaim{\noindent Proposition 2.2.}
For all $w\in H^1(\Omega\setminus\overline D)$ with $w=0$ on $\partial\Omega$
$$
\int_{\Omega\setminus\overline D}\vert w\vert^2 dy\le C_0(\Omega\setminus\overline D)^2
\int_{\Omega\setminus\overline D}\vert\nabla w\vert^2 dy
$$
where $C_0(\Omega\setminus\overline D)$ is a positive constant independent of $w$.
\endproclaim

\noindent
Another is given in the following.

\proclaim{\noindent Proposition 2.3.} Let $U$ be a bounded Lipschitz domain
of $\Bbb R^m$. For any $v\in H^1(U)$ and measurable $A\subset U$ with $\vert A\vert>0$ we have
$$
\displaystyle
\int_U\vert v-v_A\vert^2 dy\le C(U)^2(1+\frac{\vert U\vert^{1/2}}{\vert A\vert^{1/2}})^2\int_U\vert\nabla v\vert^2 dy
$$
where $C(U)$ is a positive constant independent of $v$ and $A$;
$$
\displaystyle
v_A=\frac{1}{\vert A\vert}\int_A vdy.
$$
\endproclaim

\noindent
We make use of the property that $\displaystyle C(U)^2(1+\frac{\vert U\vert^{1/2}}{\vert A\vert^{1/2}})^2$
depends on $\vert A\vert$ continuously for each fixed $U$.

\noindent
The following lemma is a special case of
Theorem 1.5.1.10 in page 41 of \cite{G}.

\proclaim{\noindent Lemma 2.4.} Let $W$ be a bounded open subset
of $\Bbb R^m$ with a Lipschiz boundary $\Gamma$. Then there exists
a positive constant $K(W)$ such that
$$\displaystyle
\int_{\Gamma}\vert u\vert^2dS
\le K(W)(\epsilon\int_{W}\vert\nabla u\vert^2dy+\epsilon^{-1}\int_W\vert u\vert^2 dy)
$$
for all $u\in H^1(W)$ and $\epsilon\in]0,\,1[$.

\endproclaim

Now we can state another side of the probe method.

\proclaim{\noindent Theorem B.}
Let $k\ge 0$ and $L>0$ satisfy
$$
\displaystyle
\Vert \text{Re}\,\lambda\Vert_{L^{\infty}(\partial D)}
+\Vert\text{Im}\,\lambda\Vert_{L^{\infty}(\partial D)}\le L;
\tag {2.6}
$$

\noindent
for some $\epsilon\in]0,\,1[$
$$\displaystyle
2K(\Omega\setminus\overline D)L\epsilon
+(k^2+2K(\Omega\setminus\overline D)L\epsilon^{-1})C_0(\Omega\setminus\overline D)^2\le 1
\tag {2.7}
$$
and
$$\displaystyle
\min_j\{1-2K(D)L\epsilon-2(k^2+2K(D)L\epsilon^{-1})C(D_j)^2(1+1)^2\}>0.
\tag {2.8}
$$
Let $x\in\Omega$ and $\sigma\in N_x$.
If $x\in\Omega\setminus\overline D$ and $\sigma$ satisfies $\sigma(]0,1])\cap D\not=\emptyset$
or $x\in\overline D$, then for any needle sequence $\xi=\{v_n\}$ for $(x,\sigma)$ we have
$\lim_{n\longrightarrow\infty}I(x,\sigma,\xi)_n=\infty$.

\endproclaim

\noindent
Note that both $\partial\Omega$ and $\partial D$ are Lipschitz;
$\lambda\in L^{\infty}(\partial D)$.
(2.7) and (2.8) mean that both $k$ and $L$ are small.

\noindent
Remarks are in order.

$\bullet$  This theorem does not cover the case when
$x\in\Omega\setminus\overline D$ and $\sigma$ satisfies both
$\sigma(]0,1])\cap D=\emptyset$ and $\sigma(]0,1])\cap\overline
D\not=\emptyset$. However, this is quite an exceptional case.

$\bullet$ At the present time we do not know how to drop the conditions (2.7) and (2.8).
This remains open.

Here we make a definition.
Let $\mbox{\boldmath $b$}$ be a nonzero vector in $\Bbb R^m$.
Given $x\in\Bbb R^m$, $\rho>0$ and $\theta\in]0,\pi[$
the set
$$\displaystyle
V=\{y\in\Bbb R^m\,\vert\,\vert y-x\vert<\rho\,\,\text{and}\,\,
(y-x)\cdot\mbox{\boldmath $b$}>\vert y-x\vert\vert\mbox{\boldmath $b$}\vert
\cos(\theta/2)\}
$$
is called a finite cone of height $\rho$, axis direction $\mbox{\boldmath $b$}$ and
aperture angle $\theta$ with vertex at $x$.

\noindent
The two lemmas given below are the core of side B of the probe method and tell us that any
needle sequence for any needle blows up on the needle.
See \cite{I6} for the proof.

\proclaim{\noindent Lemma 2.5.}
Let $x\in\Omega$ be an arbitrary point and $\sigma$ be a needle with tip at $x$.
Let $\xi=\{v_n\}$ be an arbitrary needle sequence for $(x,\sigma)$.
Then, for any finite cone $V$ with vertex at $x$ we have
$$\displaystyle
\lim_{n\longrightarrow\infty}\int_{V\cap\Omega}\vert\nabla v_n(y)\vert^2dy=\infty.
$$
\endproclaim

\proclaim{\noindent Lemma 2.6.}
Let $x\in\Omega$ be an arbitrary point and $\sigma$ be a needle with tip at $x$.
Let $\xi=\{v_n\}$ be an arbitrary needle sequence for $(x,\sigma)$.
Then for any point $z\in\sigma(]0,\,1[)$ and open ball $B$ centered at $z$ we have
$$\displaystyle
\lim_{n\longrightarrow\infty}\int_{B\cap\Omega}\vert\nabla v_n(y)\vert^2dy=\infty.
$$
\endproclaim

{\bf\noindent 2.3.  Remark}

\noindent
Finally we point out that as a corollary of (A.1) and Theorem B
we obtain a characterization of the obstacle by using the indicator sequence.

\proclaim{\noindent Corollary C.}  Under the same assumptions as those of Theorem B
we have: a point $x\in\Omega$ belongs to $\Omega\setminus\overline D$ if and only if
there exist a needle $\sigma$ with tip at $x$ and needle sequence
$\xi$ for $(x,\sigma)$ such that
the sequence $\{I(x,\sigma,\xi)_n\}$ is bounded from above.

\endproclaim

For applying side B of the probe method to non volumetric discontinuity (crack)
one needs an idea of making use of the reflected solution only to show
the blowup of the indicator sequence.  We have already developed
the idea in \cite{I7}.  In Section 4 we present an application of the
idea to the study of the reflected solutions.

\section{Proof of Theorems A and B}

{\bf\noindent 3.1.  Proof of Theorem B}

\noindent
Applying Lemma 2.4 for $W=\Omega\setminus\overline D$ or $D$ to the right hand side
of (2.3), we have
$$\begin{array}{c}
\displaystyle
\text{Re}\,<(\Lambda_0-\Lambda_D)f,\overline f>
\ge (1-2K(\Omega\setminus\overline D)L\epsilon)
\int_{\Omega\setminus\overline D}\vert\nabla(u-v)\vert^2 dy\\
\\
\displaystyle
-(k^2+2K(\Omega\setminus\overline D)L\epsilon^{-1})
\int_{\Omega\setminus\overline D}\vert u-v\vert^2dy\\
\\
\displaystyle
+(1-2K(D)L\epsilon)\int_D\vert\nabla v\vert^2 dy
-(k^2+2K(D)L\epsilon^{-1})\int_D\vert v\vert^2dy.
\end{array}
\tag {3.1}
$$

\noindent
Applying Proposition 2.2 to the second term of (3.1), we have
$$\begin{array}{c}
\displaystyle
\text{Re}\,<(\Lambda_0-\Lambda_D)f,\overline f>\\
\\
\displaystyle
\ge \{1-2K(\Omega\setminus\overline D)L\epsilon
-(k^2+2K(\Omega\setminus\overline D)L\epsilon^{-1})C_0(\Omega\setminus\overline D)^2\}
\int_{\Omega\setminus\overline D}\vert\nabla(u-v)\vert^2 dy\\
\\
\displaystyle
+(1-2K(D)L\epsilon)\int_D\vert\nabla v\vert^2 dy
-(k^2+2K(D)L\epsilon^{-1})\int_D\vert v\vert^2dy.
\end{array}
\tag {3.2}
$$
On the other hand, from Proposition 2.3 we have
$$\begin{array}{c}
\displaystyle
\int_{D}\vert v\vert^2 dy=\sum_{j}\int_{D_j}\vert v\vert^2dy\\
\\
\displaystyle
\le\sum_{j} 2\int_{D_j}\vert v-v_{A_j}\vert^2 dy
+2\int_{D_j}\vert v_{A_j}\vert^2 dy\\
\\
\displaystyle
\le\sum_{j}2C(D_j)^2(1+\frac{\vert D_j\vert^{1/2}}{\vert A_j\vert^{1/2}})^2
\int_{D_j}\vert\nabla v\vert^2 dy
+\sum_{j}2\vert D_j\vert\vert v_{A_j}\vert^2
\end{array}
$$
where $A_j\subset D_j$ and satisfy $\vert A_j\vert>0$.

\noindent Using this inequality, (2.6) and (2.7), from (3.2) we
have the basic inequality
$$\begin{array}{c}
\displaystyle
\text{Re}\,<(\Lambda_0-\Lambda_{D})f,\overline f>\\
\\
\displaystyle
\ge \sum_{j}\{1-2K(D)L\epsilon-2(k^2+2K(D)L\epsilon^{-1})C(D_j)^2(1+\frac{\vert D_j\vert^{1/2}}
{\vert A_j\vert^{1/2}})^2\}\int_{D_j}\vert\nabla v\vert^2 dy\\
\\
\displaystyle
-2(k^2+2K(D)L\epsilon^{-1})\vert D\vert\sum_{j}\vert v_{A_j}\vert^2.
\end{array}
\tag {3.3}
$$
Hereafter we proceed along the same line as \cite{I6} and completes the proof.
However, for reader's convenience we present the detail.
Choose a sequence $\{K_l\}$ of compact sets of $\Bbb R^m$ in such a way that
$K_l\subset\Omega\setminus\sigma(]0,1])$; $\overline K_l\subset K_{l+1}$ for $l=1,\cdots$;
$\Omega\setminus\sigma(]0,1])=\cup_{l=1}^{\infty}K_l$.
Then $\vert K_l\cap D_j\vert\longrightarrow\vert D_j\setminus\sigma(]0,1])\vert
=\vert D_j\vert$ as $l\longrightarrow\infty$ uniformly with $j=1,\cdots, N$.
Thus one can take a large $l_0$ in such a way that
the set $A_j\equiv K_{l_0}\cap D_j$ satisfies
$$\begin{array}{c}
\displaystyle
\max_{j}\{2(k^2+2K(D)L\epsilon^{-1})(C(D_j)^2(1+\frac{\vert D_j\vert^{1/2}}{\vert A_j\vert^{1/2}})^2
-C(D_j)^2(1+1)^2)\}\\
\\
\displaystyle
<\min_{j}\{1-2K(D)L\epsilon-2(k^2+2K(D)L\epsilon^{-1})C(D_j)^2(1+1)^2\}.
\end{array}
$$
Note that this right hand side is positive because of (2.8). We
know that the sequences $\{(v_n)_{A_j}\}$ for each $j=1,\cdots, N$
are always convergent since $\overline
A_j\subset\Omega\setminus\sigma(]0,1])$. From (3.3) with $f=v_n\vert_{\partial\Omega}$ we have
$$
\displaystyle
I(x,\sigma,\xi)_n\ge
NC\int_{D}\vert\nabla v_n\vert^2 dy
-2(k^2+2K(D)L\epsilon^{-1})\vert D\vert\sum_{j}\vert (v_n)_{A_j}\vert^2
$$
where
$$
\begin{array}{c}
\displaystyle
C=\min_{j}\{1-2K(D)L\epsilon-2(k^2+2K(D)L\epsilon^{-1})C(D_j)^2(1+1)^2\}\\
\\
\displaystyle
-\max_{j}\{2(k^2+2K(D)L\epsilon^{-1})(C(D_j)^2(1+\frac{\vert D_j\vert^{1/2}}{\vert A_j\vert^{1/2}})^2-
C(D_j)^2(1+1)^2)\}>0.
\end{array}
$$
Then the blowup of $I(x,\sigma,\xi)_n$ comes from the blowup
of the sequence
$$
\displaystyle
\int_D\vert\nabla v_n\vert^2 dy.
\tag {3.4}
$$
If $x\in D$, then the blowup of the sequence given by (3.4) is
a direct consequence of Lemma 2.5.  If $x\in\partial D$, then the
exists a finite cone $V$ at vertex at $x$ such that $V\subset D$.
This is because of the Lipshitz regularity of $\partial D$. Then
Lemma 2.5 gives the blowup of the sequence.  Now consider the
case when $x\in\Omega\setminus\overline D$ and $\sigma$ satisfies
$\sigma(]0,1])\cap D\not=\emptyset$.  Then, there exists a point $z$
on $\sigma(]0,1[)\cap D$.  Choose an open ball centered at $z$ in
such a way that $B\subset D$.  Then from Lemma 2.6 we see the
blowup of the sequence given by (3.4).

\noindent
$\Box$

In the following two subsections we always assume that $\lambda\in C^1(\partial D)$
by the reason described in subsection 2.1.

\noindent
In the following let $L>0$ satisfy (2.6).

{\bf\noindent 3.2.  First Alternative Proof of (A.3)}

\noindent
The proof is based on (3.1).

\noindent
Choose $\epsilon\in]0,\,1[$ of (3.1) in such a way that
$$
2K(\Omega\setminus\overline D)L\epsilon\le 1
$$
and
$$
2K(D)L\epsilon<1.
$$
Then (3.1) gives
$$\begin{array}{c}
\displaystyle
\text{Re}\,<(\Lambda_0-\Lambda_D)f,\overline f>
\ge -(k^2+2K(\Omega\setminus\overline D)L\epsilon^{-1})
\int_{\Omega\setminus\overline D}\vert u-v\vert^2dy\\
\\
\displaystyle
+(1-2K(D)L\epsilon)\int_D\vert\nabla v\vert^2 dy
-(k^2+2K(D)L\epsilon^{-1})\int_D\vert v\vert^2dy.
\end{array}
\tag {3.5}
$$
Now let $v=v_n$ and $x\in\Omega\setminus\overline D$.  Then from (3.5) and (A.1),
one has
$$\begin{array}{c}
\displaystyle
I(x)
\ge -(k^2+2K(\Omega\setminus\overline D)L\epsilon^{-1})
\int_{\Omega\setminus\overline D}\vert w_x\vert^2dy\\
\\
\displaystyle
+(1-2K(D)L\epsilon)\int_D\vert\nabla G_k(y-x)\vert^2 dy
-(k^2+2K(D)L\epsilon^{-1})\int_D\vert G_k(y-x)\vert^2dy.
\end{array}
\tag {3.6}
$$
The boundedness of $L^2$-norm of $w_x$ over $\Omega\setminus\overline D$
has been established as described in subsection 2.1.
From the boundedness
of the last term of (3.6) and the blowup of the middle term
we obtain the desired conclusion.

{\bf\noindent 3.3.  Second Alternative Proof of (A.3)}

\noindent
This proof does not make use of (3.1).  The key point is a weaker version of (2.5).
We prove that

\proclaim{\noindent Lemma 3.1.}
There exists a positive constant $C$ such that, for all $y_0\in \overline\Omega\setminus D$
$$\begin{array}{c}
\displaystyle
\int_{\partial D}\vert w\vert^2 dS\\
\\
\displaystyle
\le C\Vert \nabla w\Vert_{L^2(\Omega\setminus\overline D)}
\{\int_{\partial D}\vert y-y_0\vert^{1/2}
\vert\frac{\partial v}{\partial\nu}\vert dS(y)
+k^2\vert\int_{D}vdy\vert
+L\int_{\partial D}\vert v \vert dS(y)\}
\end{array}
\tag {3.7}
$$
where $v\in H^1(\Omega)$ is a solution of the equation $\triangle v+k^2v=0$ in $\Omega$;
the function $w\in H^1(\Omega\setminus\overline D)$ is the weak solution of the problem
$$\begin{array}{c}
\displaystyle
\triangle w+k^2w=0\,\,\text{in}\,\Omega\setminus\overline D,\\
\\
\displaystyle
\frac{\partial w}{\partial\nu}+\lambda w=-(\frac{\partial v}{\partial\nu}+\lambda v)
\,\,\text{on}\,\partial D,\\
\\
\displaystyle
w=0\,\,\text{on}\,\partial\Omega
\end{array}
\tag {3.8}
$$
in the following sense: the trace of $w$ onto $\partial\Omega$ vanishes and,
for all $\varphi\in H^1(\Omega\setminus\overline D)$ with $\varphi=0$ on $\partial\Omega$
$$
\displaystyle
\int_{\partial D}\lambda w \varphi dS
-\int_{\Omega\setminus\overline D}\nabla w\cdot\nabla\varphi dy
+k^2\int_{\Omega\setminus\overline D}w\varphi dy
=-\int_{\partial D}(\frac{\partial v}{\partial\nu}+\lambda v)\varphi dS.
\tag {3.9}
$$
\endproclaim

{\it\noindent Proof.}
Let $p\in H^1(\Omega\setminus\overline D)$ be the weak solution of the problem:
$$\begin{array}{c}
\triangle p+k^2 p=0\,\,\text{in}\,\Omega\setminus\overline D,\\
\\
\displaystyle
\frac{\partial p}{\partial\nu}+\lambda p=-\overline w\,\,\text{on}\,\partial D,\\
\\
\displaystyle
p=0\,\,\text{on}\,\partial\Omega.
\end{array}
\tag {3.10}
$$
It is easy to see that we have
$$\begin{array}{c}
\displaystyle
\Vert p\Vert_{H^1(\Omega\setminus\overline D)}
\le C_3\Vert w\Vert_{H^{-1/2}(\partial D)}\\
\\
\displaystyle
\le C_3\Vert w\Vert_{L^2(\partial D)}.
\end{array}
$$
Since $\lambda\in C^1(\partial D)$, we have $\partial p/\partial\nu\in H^{1/2}(\partial D)$.
Thus a standard regularity result for the Laplacian
yields $p\in H^2(\Omega\setminus\overline D)$
and the estimate
$$\displaystyle
\Vert p\Vert_{H^2(\Omega\setminus\overline D)}
\le C_4\Vert w\Vert_{H^1(\Omega\setminus\overline D)}.
$$
Then the Poincar\'e inequality gives
$$\displaystyle
\Vert p\Vert_{H^2(\Omega\setminus\overline D)}\le
C_5\Vert\nabla w\Vert_{L^2(\Omega\setminus\overline D)}.
\tag {3.11}
$$
Using (3.9) and (3.10), one can write
$$\begin{array}{c}
\displaystyle
\int_{\partial D}\vert w\vert^2 dS
=\int_{\partial D}\overline w wdS\\
\\
\displaystyle
=-\int_{\partial D}(\frac{\partial p}{\partial\nu}+\lambda p)wdS\\
\\
\displaystyle
=\int_{\Omega\setminus\overline D}(\nabla p\cdot\nabla w-k^2 pw)dy
-\int_{\partial D}\lambda p wdS\\
\\
\displaystyle
=-\int_{\partial D}p(\frac{\partial w}{\partial\nu}+\lambda w)dS\\
\\
\displaystyle
=\int_{\partial D}p(\frac{\partial v}{\partial\nu}+\lambda v)dS.
\end{array}
\tag {3.12}
$$
The Sobolev imbedding theorem ensures that $p$ can be identified with a
uniformly H\"older continuous function with exponent $1/2$ on $\overline\Omega\setminus D$
and the estimates
$$\begin{array}{c}
\displaystyle
\vert p(y)-p(y_0)\vert\le C_6\vert y-y_0\vert^{1/2}\Vert p\Vert_{H^2(\Omega\setminus\overline D)}\\
\\
\displaystyle
\vert p(y)\vert\le C_6\Vert p\Vert_{H^2(\Omega\setminus\overline D)}
\end{array}
\tag {3.13}
$$
where $C_6>0$ is independent of $p$, $y$ and $y_0$, are valid.
A combination of (3.12) and (3.13) yields
$$
\begin{array}{c}
\displaystyle
\int_{\partial D}\vert w\vert^2 dS
\le\vert\int_{\partial D}(p(y)-p(y_0))
\frac{\partial v}{\partial\nu}dS(y)\vert\\
\\
\displaystyle
+\vert p(y_0)\int_{\partial D}\frac{\partial v}{\partial\nu}dS(y)\vert
+\vert\int_{\partial D}p\lambda v dS(y)\vert\\
\\
\displaystyle
\le C_6\Vert p\Vert_{H^2(\Omega\setminus\overline D)}
\{\int_{\partial D}\vert y-y_0\vert^{1/2}
\vert\frac{\partial v}{\partial\nu}\vert dS(y)
+k^2\vert\int_{D}vdy\vert
+L\int_{\partial D}\vert v \vert dS(y)\}.
\end{array}
$$
From this and (3.11) we see that (3.7) is valid for $C=C_6C_5$.

\noindent
$\Box$

\noindent
We call the function $w$ in Lemma 3.1 the {\it reflected solution} of $v$ by $D$.

\noindent
Let $v=v_n$ and $w_n$ is the weak solution of (3.8).
From (2.3) we have
$$\begin{array}{c}
\displaystyle
I(x,\sigma,\xi)_n
\ge
\int_{\Omega\setminus\overline D}\vert\nabla w_n\vert^2dy
-k^2\int_{\Omega\setminus\overline D}\vert w_n\vert^2 dy\\
\\
\displaystyle
+\int_{D}\vert\nabla v_n\vert^2 dy-k^2\int_{D}\vert v_n\vert^2 dy
-2L\int_{\partial D}(\vert v_n\vert^2+\vert w_n\vert^2)dS\\
\\
\displaystyle
\ge (1-\epsilon CL)\int_{\Omega\setminus\overline D}\vert\nabla w_n\vert^2 dy
-k^2\int_{\Omega\setminus\overline D}\vert w_n\vert^2 dy\\
\\
\displaystyle
+\int_{D}\vert\nabla v_n\vert^2 dy-k^2\int_{D}\vert v_n\vert^2 dy
-2L\int_{\partial D}\vert v_n\vert^2 dS\\
\\
\displaystyle
-\epsilon^{-1}CL
\{\int_{\partial D}\vert y-y_0\vert^{1/2}
\vert\frac{\partial v_n}{\partial\nu}\vert dS(y)
+k^2\vert\int_{D}v_ndy\vert
+L\int_{\partial D}\vert v_n\vert dS(y)\}^2
\end{array}
\tag {3.14}
$$
where $\epsilon>0$ is arbitrary.
Let $x\in\Omega\setminus\overline D$ and take $y_0=x$ in Lemma 3.1.
Choose $\epsilon$ in such a way that $1\ge\epsilon CL$
and letting $n\longrightarrow\infty$, from (3.14) and (A.1) we obtain
$$\begin{array}{c}
\displaystyle
I(x)
\ge
-k^2\int_{\Omega\setminus\overline D}\vert w_x\vert^2 dy\\
\\
\displaystyle
+\int_{D}\vert\nabla G_k(y-x)\vert^2 dy-k^2\int_{D}\vert G_k(y-x)\vert^2 dy
\\
\\
\displaystyle
-2L\int_{\partial D}\vert G_k(y-x)\vert^2 dS
-\epsilon^{-1}CL
\{\int_{\partial D}\vert y-x\vert^{1/2}
\vert\frac{\partial}{\partial\nu}G_k(y-x)\vert dS(y)\\
\\
\displaystyle
+k^2\vert\int_{D}G_k(y-x)dy\vert
+L\int_{\partial D}\vert G_k(y-x)\vert dS(y)\}^2.
\end{array}
$$
Now one can easily see that $\lim_{x\longrightarrow a}I(x)=\infty$.

\section{Blowup of Reflected Solution}

In this section we study the blowup property
of the energy of the reflected solutions.  This is an application
of the idea developed in \cite{I7}.  The starting point is the following simple
estimate.

\proclaim{\noindent Lemma 4.1.}
Let $v\in H^1(\Omega)$ be a solution of the equation $\triangle v+k^2v=0$ in $\Omega$.
Let $w\in H^1(\Omega\setminus\overline D)$ be the reflected solution of $v$
by $D$.  This means that the function $w$ is the weak solution of (3.8).
Then we have
$$\begin{array}{c}
\displaystyle
\frac{\displaystyle\int_D\vert\nabla v\vert^2dy-k^2\int_D\vert v\vert^2dy+\int_{\partial D}
(\text{Re}\,\lambda)\vert v\vert^2dS}
{\displaystyle \Vert v\Vert_{H^1(D)}}\\
\\
\displaystyle
\le C\Vert\nabla w\Vert_{L^2(\Omega\setminus\overline D)}
\end{array}
\tag {4.1}
$$
where $C$ is a positive constant and independent of $v$ and $w$.
\endproclaim

{\it\noindent Proof.}
From the trace theorem we know that there exists $p\in H^1(\Omega\setminus\overline D)$
such that
$$\begin{array}{c}
\displaystyle
p=\overline v\,\,\text{on}\,\partial D,\\
\\
\displaystyle
p=0\,\,\text{on}\,\partial\Omega.
\end{array}
$$
This $p$ satisfies
$$\displaystyle
\Vert p\Vert_{H^1(\Omega\setminus\overline D)}
\le C_1\Vert v\vert_{\partial D}\Vert_{H^{1/2}(\partial D)}
\tag {4.2}
$$
where $C_1=C_1(\Omega\setminus\overline D)>0$ and independent of $v$ and $k^2$.
Integration by parts gives
$$\begin{array}{c}
\displaystyle
\int_D\vert\nabla v\vert^2dy-k^2\int_D\vert v\vert^2dy
=\int_{\partial D}\frac{\partial v}{\partial\nu}\overline v dS\\
\\
\displaystyle
=\int_{\partial D}(\frac{\partial v}{\partial\nu}+\lambda v)\overline v dS
-\int_{\partial D}\lambda\vert v\vert^2dS.
\end{array}
\tag {4.3}
$$
Then a combination of (3.9) for $\varphi=p$ and (4.3) yields
$$\begin{array}{c}
\displaystyle
\int_D\vert\nabla v\vert^2dy-k^2\int_D\vert v\vert^2dy+\int_{\partial D}\lambda\vert v\vert^2dS\\
\\
\displaystyle
=\int_{\Omega\setminus\overline D}\nabla w\cdot\nabla p-k^2\int_{\Omega\setminus\overline D}wpdy
-\int_{\partial D}\lambda w\overline vdS.
\end{array}
\tag {4.4}
$$

\noindent
Let $C_2=C_2(D)>0$ satisfy, for all $\Psi\in H^1(D)$
$$\displaystyle
\Vert\Psi\vert_{\partial D}\Vert_{H^{1/2}(\partial D)}\le C_2\Vert \Psi\Vert_{H^1(D)}.
\tag {4.5}
$$
Taking the real part of (4.4), from (4.5) we have
$$\begin{array}{c}
\displaystyle
\int_D\vert\nabla v\vert^2 dy-k^2\int_D\vert v\vert^2 dy+\int_{\partial D}(\text{Re}\,\lambda)\vert v\vert^2 dS\\
\\
\displaystyle
\le \Vert\nabla w\Vert_{L^2(\Omega\setminus\overline D)}\Vert\nabla p\Vert_{L^2(\Omega\setminus\overline D)}
+k^2\Vert w\Vert_{L^2(\Omega\setminus\overline D)}\Vert p\Vert_{L^2(\Omega\setminus\overline D)}
+L\Vert w\vert_{\partial D}\Vert_{L^2(\partial D)}\Vert v\vert_{\partial D}\Vert_{L^2(\partial D)}\\
\\
\displaystyle
\le C_1(\Vert\nabla w\Vert_{L^2(\Omega\setminus\overline D)}+k^2\Vert w\Vert_{L^2(\Omega\setminus\overline D)})
\Vert v\vert_{\partial D}\Vert_{H^{1/2}(\partial D)}
+L\Vert w\vert_{\partial D}\Vert_{L^2(\partial D)}\Vert v\vert_{\partial D}\Vert_{L^2(\partial D)}\\
\\
\displaystyle
\le C_1C_2(\Vert\nabla w\Vert_{L^2(\Omega\setminus\overline D)}
+k^2\Vert w\Vert_{L^2(\Omega\setminus\overline D)}+\frac{L}{C_1}\Vert w\vert_{\partial D}\Vert_{L^2(\partial D)})
\Vert v\Vert_{H^1(D)}
\end{array}
\tag {4.6}
$$
Now applying Proposition 2.1, the trace theorem in the domain $\Omega\setminus\overline D$ to
the factor of (4.6) involving $w$,
we obtain (4.1).

$\Box$

\noindent
The point is: everything has been done in the context of the weak solution.  We do not
make use of higher regularity.

Let $L>0$ satisfy (2.6).
Then from Lemma 2.4 we have, for all $\epsilon\in\,]0,\,1[$
$$\begin{array}{c}
\displaystyle
\int_D\vert\nabla v\vert^2dy-k^2\int_D\vert v\vert^2dy+\int_{\partial D}
(\text{Re}\,\lambda)\vert v\vert^2dS\\
\\
\displaystyle
\ge (1-K(D)L\epsilon)\int_D\vert\nabla v\vert^2 dy
-(k^2+K(D)L\epsilon^{-1})\int_D\vert v\vert^2dy.
\end{array}
\tag {4.7}
$$

{\bf\noindent 4.1.  Blowup of Energy of Reflected Solution by $D$}

\noindent
Given $x\in\Omega\setminus\overline D$ let $w_x$ be the reflected solution by $D$.
One can choose a needle $\sigma$ with tip at $x$ in such a way that
$\sigma(]0,\,1])\cap\overline D=\emptyset$.  Then any needle sequence
$\xi=\{v_n\}$ for $(x,\sigma)$ is convergent on $D$.  Then the sequence $\{w_n\}$ of reflected solutions
$w_n$ of $v_n$ by $D$ converges to $w_x$.  Thus from (4.1) and (4.7) we obtain the estimate
$$\begin{array}{c}
\displaystyle
C\Vert\nabla w_x\Vert_{L^2(\Omega\setminus\overline D)}\\
\\
\displaystyle
\ge
\frac{\displaystyle(1-K(D)L\epsilon)\int_D\vert\nabla G_k(y-x)\vert^2dy-
(k^2+K(D)L\epsilon^{-1})\int_D\vert G_k(y-x)\vert^2dy}
{\displaystyle(\int_D\vert\nabla G_k(y-x)\vert^2 dy+\int_D\vert G_k(y-x)\vert^2dy)^{1/2}}.
\end{array}
$$
Thus choosing $\epsilon$ in such a way that $1-LK(D)\epsilon>0$, one concludes that,
for all $a\in\partial D$
$$
\lim_{x\longrightarrow a}\int_{\Omega\setminus\overline D}\vert\nabla w_x\vert^2 dy=\infty.
$$

{\bf\noindent 4.2.  Blowup of Energy of Reflected Solution of $v_n$ by $D$}

\noindent
The following theorem gives a characterization of the blowup
of the energy of $v_n$ on $D$ by the blowup of the energy
of the reflected solution of $v_n$  by $D$.

\proclaim{\noindent Theorem D.}
Let $k$, $L$ and $\epsilon\in\,]0,\,1[$ satisfy
$$\displaystyle
\min_j\{1-K(D)L\epsilon-2(k^2+K(D)L\epsilon^{-1})C(D_j)^2(1+1)^2\}>0.
\tag {4.8}
$$
Let $x\in\Omega$ and $\sigma$ be a needle with tip at $x$. Let $\xi=\{v_n\}$ be a
needle sequence for $(x,\sigma)$. Let $\{w_n\}$
be the sequence of the reflected solutions of $v_n$ by $D$.
If
$$\displaystyle
\lim_{n\longrightarrow\infty}\int_{D}\vert\nabla v_n\vert^2dy=\infty,
$$
then
$$
\lim_{n\longrightarrow\infty}\int_{\Omega\setminus\overline D}\vert\nabla w_n\vert^2 dy=\infty.
$$
\endproclaim

\noindent
Note that (2.8) in Theorem B implies (4.8).

{\it\noindent Proof.} Applying the argument after (3.2) in the
proof of Theorem B to the right hand side of (4.7), we obtain the
estimate
$$\begin{array}{c}
\displaystyle
(1-K(D)L\epsilon)\int_D\vert\nabla v_n\vert^2 dy
-(k^2+K(D)L\epsilon^{-1})\int_D\vert v_n\vert^2dy\\
\\
\displaystyle
\ge
C'\int_D\vert\nabla v_n\vert^2 dy-C''
\end{array}
\tag {4.9}
$$
where $C'$ and $C''$ are positive constant independent of $n$.

\noindent
Here we cite a lemma in \cite{I7} that says that $\nabla v_n$ dominates $v_n$.

\proclaim{\noindent Lemma 4.2(\cite{I7}).} Let $x\in\Omega$ and
$\sigma$ be a needle with tip at $x$. Let $\xi=\{v_n\}$ be a
needle sequence for $(x,\sigma)$. If
$$\displaystyle
\lim_{n\longrightarrow\infty}\int_{D}\vert\nabla v_n\vert^2dy=\infty,
$$
then there exists a natural number $n_0$ such that the sequence
$$
\displaystyle \{\frac{\displaystyle\int_{D}\vert
v_n\vert^2dy} {\displaystyle\int_{D}\vert\nabla
v_n\vert^2dy}\}_{n\ge n_0},
$$
is bounded.
\endproclaim

Now from (4.1), (4.7), (4.9) and Lemma 4.2 yields Theorem D.
The converse of Theorem D is also true without the condition (4.8).
It is an easy consequence of the well posedness of (2.2), the trace theorem in the domain $D$ and Lemma 4.2.

\section{Enclosure Method}

In this section we present an application of the enclosure method
introduced in \cite{IE} to {\it Inverse Problem} in three-dimensions ($m=3$).

\noindent In \cite{IE} we considered the problem in the case when
$\lambda\longrightarrow\infty$. The case when $\lambda\longrightarrow 0$ was studied in
\cite{IE2}. Therein we demonstrated how to use the exponentially
growing solutions for elliptic equations constructed by
Sylvester-Uhlmann \cite{SU} for extracting the convex hull of
unknown obstacles from the Dirichlet-to-Neumann map. Here we apply
the idea to {\it Inverse Problem}.

The convex hull of $D$ is uniquely determined by the support function
$$\displaystyle
h_D(\omega)=\sup_{x\in\, D}x\cdot\omega,\,\,\omega\in S^2.
$$
The exponentially growing solutions for the Helmholtz equation
have the form
$$\displaystyle
v(x)=e^{x\cdot(\tau\omega+i\sqrt{\tau^2+k^2}\omega^{\perp})}
\tag {5.1}
$$
where $\omega^{\perp}\in S^2$ and satisfies $\omega\cdot\omega^{\perp}=0$;
$\tau>0$ is a large parameter.

Using the traces of these functions onto $\partial\Omega$,
we define the indicator function
$$\displaystyle
I_{\omega,\,\omega^{\perp}}(\tau,t)
=e^{-2\tau t}\text{Re}\,<(\Lambda_0-\Lambda_D)f,\overline f>, -\infty<t<\infty
$$
where $f=v$ on $\partial\Omega$.

Here we prove the following.

\proclaim{\noindent Theorem E.}  Assume that both $\partial\Omega$ and $\partial D$ are $C^2$
and $\lambda\in C^1(\partial D)$.  Let $\omega$ and $\partial D$ satisfy the condition
(R): the set $\{x\in\partial D\,\vert\,x\cdot\omega=h_D(\omega)\}$
consists of a single point and the Gaussian curvature of $\partial D$ does not vanish
at the point.

\noindent
Then the formula
$$\displaystyle
\lim_{\tau\longrightarrow\infty}
\frac{\log\vert I_{\omega,\,\omega^{\perp}}(\tau,0)\vert}{2\tau}=h_D(\omega),
\tag {5.2}
$$
is
valid.  Moreover we have:

if $t>h_D(\omega)$, then $\lim_{\tau\longrightarrow\infty}\vert I_{\omega,\,\omega^{\perp}}(\tau,t)\vert=0$;

if $t=h_D(\omega)$, then $\liminf_{\tau\longrightarrow\infty}\vert I_{\omega,\,\omega^{\perp}}(\tau,t)\vert>0$;

if $t<h_D(\omega)$, then $\lim_{\tau\longrightarrow\infty}\vert I_{\omega,\,\omega^{\perp}}(\tau,t)\vert=\infty$.

\endproclaim

{\it\noindent Proof.}
The proof can be done along the same line as Theorem 4.2 of \cite{IE2}.
Thus it suffices to prove that there exists a positive constant $C$ independent of $\tau>>1$
such that, for all $\tau>>1$
$$\displaystyle
\vert I_{\omega,\,\omega^{\perp}}(\tau,h_D(\omega))\vert\ge C.
\tag {5.3}
$$
This is proved as follows.
From (3.1) for $v$ given by (5.1) we obtain
$$\begin{array}{c}
\displaystyle
\vert I_{\omega,\,\omega^{\perp}}(\tau,h_D(\omega))\vert
\ge -(k^2+2K(\Omega\setminus\overline D)L\epsilon^{-1})e^{-2\tau h_D(\omega)}
\int_{\Omega\setminus\overline D}\vert w\vert^2dy\\
\\
\displaystyle
+\{(1-2K(D)L\epsilon)-\frac{k^2+2K(D)L\epsilon^{-1}}
{2\tau^2+k^2}\}e^{-2\tau h_D(\omega)}\int_D\vert\nabla v\vert^2dy
\end{array}
\tag {5.4}
$$
where $L>0$ satisfy (2.6) and $w$ is the weak solution of the problem (3.8).
We claim
$$\displaystyle
\lim_{\tau\longrightarrow\infty}
\frac{\displaystyle e^{-2\tau h_D(\omega)}\int_{\Omega\setminus\overline D}\vert w\vert^2dy}
{\displaystyle e^{-2\tau h_D(\omega)}\int_D\vert\nabla v\vert^2dy}=0.
\tag {5.5}
$$
This is proved as follows.

\noindent
Let $y_0$ be the single point in the condition (R).
Using the completely same argument as \cite{IE2}, one obtains the estimate
$$\begin{array}{c}
\displaystyle
e^{-2\tau h_D(\omega)}\int_{\Omega\setminus\overline D}\vert w\vert^2dy\\
\\
\displaystyle
\le
C'\{(\int_{\partial D}\vert y-y_0\vert^{1/2}e^{-\tau h_D(\omega)}
\vert\frac{\partial v}{\partial\nu}+\lambda v\vert dS(y))^2\\
\\
\displaystyle
+\vert k^2\int_De^{-\tau h_D(\omega)}v dy-\int_{\partial D}\lambda e^{-\tau h_D(\omega)}v dS\vert^2\}
\end{array}
\tag {5.6}
$$
where $C'$ is a positive constant independent of $\tau$.  From (17) of \cite{IE2}
we have immediately,
as $\tau\longrightarrow\infty$
$$
(\int_{\partial D}\vert y-y_0\vert^{1/2}e^{-\tau h_D(\omega)}
\vert\frac{\partial v}{\partial\nu}+\lambda v\vert dS(y))=O(\tau^{-1/2}).
\tag {5.7}
$$
Here we made use of the non vanishing of the Gaussian curvature of
$\partial D$ at $y_0$. From the Schwartz inequality we have, as
$\tau\longrightarrow\infty$
$$\displaystyle
\frac{\displaystyle\vert\int_De^{-\tau h_D(\omega)}v dy\vert^2}
{\displaystyle e^{-2\tau h_D(\omega)}\int_D\vert\nabla v\vert^2dy}
=O(\tau^{-2})
\tag {5.8}
$$
and we have already known that
$$\displaystyle
e^{-2\tau h_D(\omega)}\int_D\vert\nabla v\vert^2dy\ge C'',\,\,\tau>>1
\tag {5.9}
$$
where $C''$ is a positive constant independent of $\tau$.  Here we made use of $C^2$-regularity
of $\partial D$ at $y_0$.
Moreover we have
$$\begin{array}{c}
\displaystyle
\vert\int_{\partial D}\lambda e^{-\tau h_D(\omega)}v dS\vert^2
\le
\Vert\lambda\Vert_{L^{\infty}(\partial D)}^2(\int_{\partial D}e^{\tau(y\cdot\omega-h_D(\omega))}dS(y))^2\\
\\
\displaystyle
\le
\Vert\lambda\Vert_{L^{\infty}(\partial D)}^2\vert\partial D\vert
\int_{\partial D}e^{2\tau(y\cdot\omega-h_D(\omega))}dS(y).
\end{array}
\tag {5.10}
$$
In the proof of Lemma 4.2 of \cite{IE2} we have already proved that
$$\displaystyle
\lim_{\tau\longrightarrow\infty}
\frac{\displaystyle\int_{\partial D}e^{2\tau(y\cdot\omega-h_D(\omega)}dS(y)}
{\displaystyle e^{-2\tau h_D(\omega)}\int_D\vert\nabla v\vert^2dy}=0.
\tag {5.11}
$$
Now from (5.6) to (5.11) one obtains (5.5).
A combination of (5.4) for sufficiently small $\epsilon$, (5.5) and (5.9) gives (5.3).

\noindent
$\Box$

$$\quad$$

\centerline{{\bf Acknowledgement}}

This research was partially supported by Grant-in-Aid for Scientific
Research (C)(2) (No.  15540154) of Japan  Society for the Promotion of Science.

$$\quad$$

\vskip1cm
\noindent
e-mail address

ikehata@math.sci.gunma-u.ac.jp

\begin{thebibliography}{99}


\bibitem{CLN} Cheng, J., Liu, J. and Nakamura, G.,
           \newblock  Recovery of the shape of an obstacle and the boundary impedance from the far-field pattern,
           \newblock  J. Math. Kyoto Univ., 43(2003), no. 1, 165-186.




\bibitem{CLN2} Cheng, J., Liu, J. J. and Nakamura, G.,
           \newblock The numerical realization of the probe method for the inverse scattering problems from the near
           field data,
           \newblock Inverse Problems, 21(2005), 839-855.










\bibitem{CK} Colton, D. and Kress, R.,
           \newblock Eigenvalues of the far field operator for the Helmholtz
           equation in an absorbing medium,
           \newblock SIAM J. Appl. Math., 55(1995), 1724-1735.



\bibitem{CP} Colton, D. and Piana, M.,
           \newblock Inequalities for inverse scattering problems in absorbing media,
           \newblock Inverse Problems, 17(2001), 597-605.




\bibitem{DIN} Daido, Y., Ikehata, M. and Nakamura, G.,
           \newblock Reconstruction of inclusions for the inverse boundary value problem with mixed type
           boundary condition,
           \newblock Appl. Anal., 83(2004), 109-124.









\bibitem{EP} Erhard, K. and Potthast, R.,
          \newblock A numerical study of the probe method, submitted.



\bibitem{GK} Grinberg, N. and Kirsch, A.,
           \newblock The linear sampling method in inverse obstacle scattering problem
           for impedance boundary conditions,
           \newblock J. Inv. Ill-Posed Problems, 10(2002), 171-185.



\bibitem{G} Grisvard, P.,
          \newblock  Elliptic problems in nonsmooth domains, Pitman, Boston, 1985.




\bibitem{I1} Ikehata, M.,
          \newblock Reconstruction of the shape of the inclusion by boundary measurements,
          \newblock Comm. PDE., 23(1998), 1459-1474.



\bibitem{I2} Ikehata, M.,
          \newblock Reconstruction of an obstacle from the scattering amplitude at a fixed frequency,
          \newblock Inverse Problems, 14(1998), 949-954.



\bibitem{I3} Ikehata, M.,
             \newblock Reconstruction of obstacle from boundary measurements,
             \newblock Wave Motion, 30(1999), 205-223.


\bibitem{IE} Ikehata, M.,
             \newblock Reconstruction of the support function for inclusion from boundary
             measurements,
             \newblock J. Inv. Ill-Posed Problems, 8(2000), 367-378.


\bibitem{IE2} Ikehata, M.,
             \newblock How to draw a picture of an unknown inclusion from boundary measurements.
             Two mathematical inversion algorithms,
             \newblock J. Inv. Ill-Posed Problems, 7(1999), 255-271.



\bibitem{I4} Ikehata, M.,
             \newblock Reconstruction of inclusion from boundary measurements,
             \newblock J. Inv. Ill-Posed Problems, 10(2002), 37-65.



\bibitem{I6} Ikehata, M.,
             \newblock A new formulation of the probe method and related problems,
             \newblock Inverse Problems, 21(2005), 413-426.


\bibitem{I7} Ikehata, M.,
             \newblock Inverse crack problem and probe method,
             \newblock Cubo (1) 8(2006), 29-40.




\bibitem{IN} Ikehata, M. and Nakamura, G.,
             \newblock Reconstruction formula for identifying cracks,
                          \newblock J. Elasticity,  70(2003), 59-72.



\bibitem{NUW} Nakamura, G., Uhlmann, G. and Wang, J. N.,
             \newblock Reconstruction of cracks in an inhomogeneous anisotropic elastic medium,
             \newblock  J. Math. Pures Appl., 82(2003), 1251-1276.





\bibitem{P1} Potthast, R.,
             \newblock A point-source method for inverse acoustic and electromagnetic obstacle scattering
             problems,
             \newblock IMA J. Appl. Math., 61(1998), 119-140.


\bibitem{P2} Potthast, R.,
             \newblock Stability estimates and reconstructions in inverse acoustic scattering using singular
             sources,
             \newblock J. Comput. Appl. Math., 114(2000), 247-274.



\bibitem{SS} Stanoyevitch, A. and Stegenga, D. A.,
             \newblock Equivalence of analytic and Sobolev Poincar\'e inequalities for
             planar domains,
             \newblock Pacific J. Math., 178(1997), 363-375.





\bibitem{SU} Sylvester, J. and Uhlmann, G.,
             \newblock A global uniqueness theorem for an inverse boundary value problem, Ann. Math. 125(1987),
             \newblock 153-169.




\bibitem{Z} Ziemer, W. P.,
            \newblock Weakly differentiable functions,
            \newblock Graduate texts in mathematics, 120,
            Springer, New York, 1989.




\end{thebibliography}
\end{document}